# Deep Domain Decomposition Method for Solving the Variational Inequality Problems

**Yiyang Wang, Qijia Zhou, Shengyuan Deng, Chenliang Li\***

*School of Mathematics and Computing Science, Guilin University of Electronics Technology, Center for Applied Mathematics of Guangxi (GUET), Guilin, Guangxi, China*
*\*Corresponding Author.*

**Abstract:** By integrating physics-informed neural network (PINN) techniques with domain decomposition method, a deep domain decomposition method is presented for solving elliptic variational inequality problems. Based on the Ritz variation method, the elliptic variational inequality problem is firstly reformulated as an optimization problem, and then the subproblem in each subdomain is solved by using the Ritz-PINN method, which the parameters in the network are updated by the Adam optimizer, and the residual-adaptive training by introducing a residual-adaptive dataset update strategy to gradually guide the model to learn more complex regions. Additionally, the impact of overlapping regions on the performance of the new algorithm is explored. Numerical results demonstrate the effectiveness of the proposed algorithm, the mean square error can be reached O (1.0e-07), and the number of iterations is independent of grid length h under uniform overlap conditions.

## 1. Introduction

Compared to traditional numerical optimization methods, the PINN model offers advantages such as high computational efficiency, low data requirements, and strong capabilities in handling nonlinearity and uncertainty. Currently, deep neural networks have achieved remarkable success in solving partial differential equations (PDEs). Weinan E [1] and colleagues employed deep learning techniques to handle general high-dimensional parabolic PDEs by approximating the gradient of the unknown solution using neural networks. Sirignano and Spiliopoulos [2] used deep approximations to solve high-dimensional PDEs by training neural networks to satisfy differential operators, initial conditions, and boundary conditions. M. Raissi [3] and collaborators introduced physics-based neural networks, incorporating both continuous-time and discrete-time models based on the distribution and arrangement of available data. Long [4] and colleagues used large multi-layer convolutional deep networks to learn PDEs from data.

The domain decomposition algorithm is a technique that divides complex problems into several simpler subproblems. In the integration of machine learning with domain decomposition, Kopaničáková [5] and colleagues constructed a nonlinear preconditioner using the Schwarz domain decomposition framework, hierarchically decomposing network parameters to provide more accurate solutions for underlying PDEs. [6] proposed a deep domain decomposition method for PDEs based on variational principles. [7] applied the deep domain decomposition method to solve elliptic problems, and numerical examples demonstrated that the iteration count is independent of the network architecture and decreases as the overlap size increases. [8] designed D3M with a hierarchical neural network framework for optimization problems. By decomposing the PDE system into component parts, the method independently constructs local neural networks within physical subdomains, enabling efficient neural network approximations for complex problems. [9] roposed a Fourier-feature-based deep domain decomposition method (F-D3M) for PDEs, which applying overlapping domain decomposition to simplify high-frequency modes into relatively low-frequency ones. In each local subdomain, a Multi-Fourier Feature

Network (MFFNet) was constructed, where effective boundary and interface treatments were applied to the corresponding loss functions. A key limitation of PINNs is lack of accuracy and efficiency when PINNs solves the larger domains and complex multi-scale problems, [10] proposed a finite-basis physics-informed neural network (FBPINN) with Schwarz domain decomposition method to accelerate the learning process of PINNs and to improve their accuracy. [11] formulated a general algorithmic framework for hp-variational physics-informed neural networks (hp-VPINN) by combining nonlinear approximations from shallow and deep neural networks with hp refinement through domain decomposition and projection into high-order polynomial spaces, effectively optimizing network parameters. [12] extended the optimized Schwarz domain decomposition method to unstructured grid problems, utilizing graph convolutional neural networks (GCNNs) and unsupervised learning to learn optimal modifications of subdomain interfaces. By improving the loss function, they achieved strong performance on arbitrarily large problems, with computational costs scaling linearly with problem size. The performance of the learned linear solver was compared with traditional and optimized domain decomposition algorithms for both structured and unstructured grid problems. [13] used neural network-enhanced operator compensation in the deep Ritz method to ensure accurate flux transmission across subregion interfaces, and presented a novel learning algorithm for a more generalized non-overlapping domain decomposition method under overfitting interface conditions.

For variational inequality problems, there are many efficient numerical methods [14-16]. Because the PINN model can embed constraint conditions and objective functions into the neural network, in recent years, machine learning methods based on Physics-Informed Neural Networks (PINNs) have provided new methods for solving variational inequality problems, such as neural-network model for monotone linear asymmetric variational inequalities[17], global exponential stability of neural networks[18], general projection neural network[19], novel neural network for variational inequalities with linear and nonlinear constraints [20] , and so on.

Although the deep domain decomposition method has achieved significant success in solving partial differential equation (PDE) problems, there has been limited research on combining the PINN model with domain decomposition methods to address variational inequality problems. Therefore, in this paper we a deep domain decomposition method to solve elliptic variational inequality problems.

The structure of this paper is as follows: Section 2 presents the deep domain decomposition method, Section 3 demonstrates the effectiveness of the method through numerical experiments, and Section 4 gives a summarize.

## 2 Deep Domain Decomposition Method Based on Improved PINNs

### 2.1 Variational Inequality Problems

Consider the following elliptic variational inequality problem: find $u \in V = \{v \in H_0^1(\Omega), v \geq 0 \text{ a.e. } \Omega\}$ such that:
$$v \in V \quad a(u, v - u) \geq f(v - u) \quad (1)$$
where
$$a(u,v) = \int_\Omega (\nabla u \cdot \nabla v + \alpha uv) d\Omega,$$
$$f \in L^2(\Omega) \quad (2)$$
And $\Omega$ is a convex domain.

By applying the Ritz variational method, we obtain the following equivalent minimization functional problem: Find $u \in V$ such that:
$$J(u) = \min J(v), \quad v \in V, \quad (3)$$
where
$$J(v) = \frac{1}{2} a(v,v) - f(v) \quad (4)$$

It is easy to see that $a(\cdot,\cdot)$ is a continuous, symmetric, and coercive bilinear functional on $\mathcal{V}$, and $f : \mathcal{V} \to \mathbf{R}$ is a continuous linear functional.

### 2.2 Deep Domain Decomposition Method (Deep DDM)

We divide the domain $\Omega$ into N subdomains, with the corresponding subproblems formulated as:
$$\begin{cases} a(u_n, v_n - u_n) \geq f(v_n - u_n), \text{ in } \Omega_n, n = 1,2,\dots,N \\ \qquad \mathcal{B}(u_n) = g_n, \text{ on } \partial\Omega \setminus \Gamma_n, \\ \qquad \mathcal{D}(u_n) = \mathcal{D}(u_n), \text{ on } \Gamma_n \\ \qquad u_n \geq 0 \end{cases} \quad (5)$$

Here, $\Gamma_n$ represents the artificial boundary between the subdomain $\Omega_n$ and other subdomains, while $\mathcal{D}$ denotes the artificial operator.

$$\text{Let } h_n(x; \theta_n), \quad (6)$$

1N and un represent the solution for each subdomain, and let $\theta\_n$ be the network parameters corresponding to the subproblem. Then, the optimal functional problem for each subdomain can be formulated as:

$$\theta^* = \text{argmin}_\theta \mathcal{M}_n(\theta; \mathbf{X}_n), \quad (7)$$

where

$$\mathcal{M}_n(\theta; \mathbf{X}_n) = w_1 \mathcal{M}_{\Omega_n}(\theta; \mathbf{X}_{f_n}) + w_2 \mathcal{M}_{\partial\Omega_n \setminus \Gamma_n}(\theta; \mathbf{X}_{g_n}) + w_3 \mathcal{M}_{\Gamma_n}(\theta; \mathbf{X}_{\Gamma_n}) + w_4 M^+(\theta), \quad (8)$$

$$\mathcal{M}_{\Omega_n}(\theta; \mathbf{X}_{fn}): \& \frac{1}{N_{fn}} \sum_{i=1}^{N_{fn}} \left| a(h_n(\mathbf{x}_{fn}^i; \theta), \mathbf{x}_{fn}^i) - f(\mathbf{x}_{fn}^i) \right|^2, \quad (9)$$

$$\mathcal{M}_{\partial\Omega_n \setminus \Gamma_n}(\theta; \mathbf{X}_{g_n}): \& \frac{1}{N_{g_n}} \sum_{i=1}^{N_{gn}} \left| \mathcal{B}(h_n(\mathbf{x}_{g_n}^i; \theta)) - g(\mathbf{x}_{g_n}^i) \right|^2 \quad (10)$$

$$\mathcal{M}_{\Gamma_n}(\theta; \mathbf{X}_{\Gamma_n}): \& = \frac{1}{N_{\Gamma n}} \sum_{i=1}^{N_{\Gamma n}} \left| \mathcal{D}(h_n(\mathbf{x}_{\Gamma_n}^i; \theta)) - \mathcal{D}(h_r(\mathbf{x}_{\Gamma_n}^i; \theta)) \right|^2, \quad (11)$$

$$M^+(\theta) := \frac{1}{N_n} \sum_{i=1}^{N_f} \max\{-u(\mathbf{x}_f^i; \theta), 0\}, \quad (12)$$

Her $\mathbf{X}_{f_n} := \{\mathbf{x}_{f_n}^i\}_{i=1}^{N_{fn}}, \mathbf{X}_{g_n} := \{\mathbf{x}_{g_n}^i\}_{i=1}^{N_{g_n}}$ and $\mathbf{X}_{\Gamma_n} := \{\mathbf{x}_{\Gamma_n}^i\}_{i=1}^{N_{\Gamma n}}$ represent the interior points, local boundary points, and interface points of the nth subdomain, respectively. $\mathbf{X}_n = \{\mathbf{X}_{f_n}, \mathbf{X}_{g_n}, \mathbf{X}_{\Gamma_n}\}$ denotes the training dataset used for each subdomain. We denote $W_s := \mathcal{D}(h_r(\mathbf{x}_{\Gamma_n}; \theta))$ as the information that would be transported to the objective subproblem labelled by s from the Neighboring subproblems labelled by r.

We construct a corresponding Physics-Informed Neural Network (PINN) within each subdomain for training. During the training process, firstly we need to minimize the loss function $\mathcal{M}_n(\theta; \mathbf{X}_n)$ through updating the parameters in the network. We choose to use the Adam optimizer (Adaptive Moment Estimation). Secondly, we introduce a residual-adaptive dataset update strategy. In the residual-adaptive training process, the model dynamically focuses on areas with larger residuals (prediction errors). Additionally, as training progresses, points in the dataset can be dynamically adjusted, gradually guiding the model to learn more complex regions. With this strategy, the model can more efficiently capture the important features in the problem, thereby improving the overall training effectiveness and solution accuracy.

In summary, we present the following Algorithm 1.

**Algorithm 1:** Improved PINNs Method

Initialize:
1. define the Physics-Informed Neural Network(PINN) structure
2. Divide the domain into n subdomains $\Omega_n$
2. weights $w_1, w_2, w_3$ for loss function components
3. construct loss function ←Variational Inequality with domain decomposition
4. Adam optimizer parameters $\alpha$
5. Bayesian optimization setup $w_1, w_2, w_3, w_4$
6. Randomly sample training data points $(x, y)$ from the domain $\Omega$

While ($epoch$):
1. Train PINN using the Adam optimizer and data points
   - Compute loss = $w_1 \mathcal{M}_{\Omega_n} + w_2 \mathcal{M}_{\partial\Omega_n \setminus \Gamma_n} + w_3 \mathcal{M}_{\Gamma_n} + w_4 M^+$,

Where:

**Region Interior Loss:**

$$M_{\Omega_n}(\theta) = \frac{1}{N_{f_n}} \sum_{i=1}^{N_{fn}} \left| \mathcal{L}\left(h_n(x_{f_n}^i, \theta)\right) - f(x_{f_n}^i) \right|^2$$

**Variational Inequality Constraint Loss:**

$$M^+(\theta) = \frac{1}{N_{f_n}} \sum_{i=1}^{N_{fn}} \max\{-u(x_{f_n}^i, \theta), 0\}$$

**Boundary Condition:**

$$M_{\partial\Omega_n \setminus \Gamma_n}(\theta) = \frac{1}{N_{g_n}} \sum_{i=1}^{N_{g_n}} \left| B\left(h_n(x_{g_n}^i, \theta)\right) - g(x_{g_n}^i) \right|^2$$

**Subdomain Interface Matching Loss:**

$$M_{\Gamma_n}(\theta) = \frac{1}{N_{\Gamma_n}} \sum_{i=1}^{N_{\Gamma_n}} \left| D\left(h_n(x^i_{\Gamma_n}, \theta)\right) - D\left(h_{r_n}(x^i_{\Gamma_n}, \theta)\right) \right|^2$$

- Update network parameters using Adam optimization steps

2. Optimize the loss weights $w_1, w_2, w_3, w_4$ using Bayesian Optimization
    - Fit surrogate model to observed losses
    - Propose new weights $w_1, w_2, w_3, w_4$ by maximizing acquisition function
    - Updata dataset for Bayesian Optimization
3. Update training dataset (under the condition that it is satisfied).
4. Evaluate stopping criteria(max epoch or loss tolerance)

end while

return: $u(x, y)$, which is an approximate solution constructed by Deep PINNs

Next, we present a domain decomposition algorithm incorporating PINNs. This algorithm divides the computational domain into multiple subdomains, where each subdomain independently trains a PINN to approximate the solution while exchanging information at interfaces.

**Algorithm 2:** Deep DDM

1. Construct $\mathbf{X}_s$;
2. Initial PINN parameters $\theta_s^0$ and interface information $W_s^0$ along $\Gamma_s$;
3. **for** $i = 1,2, \ldots$ **do**　　　　　　　　　　→Start DDM iteration
4. 　　Set $\mu_s^{i,0} := \mu_s^{i-1}$;
5. 　　**for** $j = 1,2, \ldots$ **do**　　　　　　　→Start PINNs training iteration
6. 　　　　Set $\mu_s^{i,j} := \mu_s^{i,j-1}$;
7. 　　　　Rearrange randomly training data $\{X_k^k\}_{k=1}^{m_s}$;
8. 　　　　**for** $k = 1,2, \ldots, m_s$ **do**　　→Update on minibatch
9. 　　　　　　Compule PINN loss $L_{PINN}$
10. 　　　　　Update $\mu_s^{i,j}$ via Adam
11. 　　　　**End for**
12. 　　　　**If** $\left| L_{PINN}(\mu_s^{i,j}; \mathbf{X}_s) - \mathcal{M}_s(\mu_s^{i,j-n}; \mathbf{X}_s) \right| / \left| L_{PINN}(\mu_s^{i,j}; \mathbf{X}_s) \right| < tol_L$ **then**
13. 　　　　　　BREAK;
14. 　　　　**end if**
15. 　　**end for**
16. 　　Set $\mu_s^i := \mu_s^{i,j}$;
17. 　　Interchange the interface information: $W_s^i \leftarrow \mu_s(X_{\Gamma_s})$;
18. 　　**If** $\|W_s^i - W_s^{i-1}\| / \| W_s^i \| < tol_\Gamma, \forall s$　**then**
19. 　　　　STOP;
20. 　　**If** $\|h_s(X_f; \theta_s^i) - h_s(X_f; \theta_s^{i-1})\| / \| h_s(X_f; \theta_s^i) \| < tol_\Omega, \forall s$ **then**
21. 　　　　STOP;
22. 　　**end if**
23. 　　　　　　　　　　　　　　　　　　　**end for**

## 3. Numerical Examples

In this section, the network architecture is designed with 32 units in each layer. We chose ReLU as the activation function and used the Adam optimizer with stochastic gradient descent to update the network parameters in PINNs.

Besides, the iteration count for each subdomain is set to 1k. The learning rate is adjusted based on different conditions, with a default setting of 0.001. The initial training data is randomly generated following a normal distribution, and the dataset is updated after each iteration. The initial weights are set to $10^4$.

Additionally, the computational setup includes a CPU with 192 cores running at a frequency of 1130.5 MHz, along with an RTX 4090D (24GB) graphics processor, which facilitates the high-performance training and optimization processes.

**Example 1** Consider the following inequality problem:

$$\begin{cases} -u_{xx} - u_{yy} \geq 2\pi^2 sin\pi x sin\pi y, (x,y) \in \Omega = (-1,1) \times (-1,1) \\ u \geq 0, \\ u(-u_{xx} - u_{yy} - 2\pi^2 sin\pi x sin\pi y) = 0. \\ u = 0, \ on \ \partial\Omega \end{cases} \quad (12)$$

By applying Green's formula transformation and introducing the functional J (u), the variational inequality problem is converted minimization problem.

$$\min J(u) = \frac{1}{2}\int_0^1 \int_0^1 [(u_x)^2 + (u_y)^2]dxdy - \int_0^1 \int_0^1 2\pi^2 \sin(\pi x)\sin(\pi y) \cdot u dx dy \quad (13)$$

Suppose h is the level of refinement, $\delta$ is the size of the overlap area. Let $\Omega_1 = (-1, \delta/2) \times (-1,1)$ and $\Omega_2 = (-\delta/2, 1) \times (-1,1)$, The termination condition is set as the difference between the two iterative solutions in sequence is less than $\varepsilon=10^{-4}$.

The numerical results of Algorithm 2 (Deep DDM) are shown in Table 1 and Table 2.

**Table 1. Iteration Number of Deep DDM**

| δ\h | 0.05 | 0.02 | 0.01 | 0.005 | 0.002 | 0.001 |
|---|---|---|---|---|---|---|
| h | 78 | 86 | 87 | 88 | 115 | 136 |
| 2h | 97 | 55 | 57 | 60 | 94 | 113 |
| 0.1 | 90 | 78 | 70 | 83 | 92 | 105 |
| 0.2 | 107 | 98 | 85 | 97 | 74 | 89 |

From Table 1, numerical analysis shows that the number of iterations increases with decreasing refinement h for the small overlap case (i.e., the h-overlap case and the 2h-overlap case), provided that the overlap domain $\delta$ is the same; In the consistent overlap case (i.e., the $\delta = 0.1$ and $\delta = 0.2$ cases), the number of iterations is independent of the degree of refinement h. The number of iterations is the same for all cases.

With the same degree of refinement h, the number of iterations is independent of the overlap domain $\delta$ for the small overlap case (i.e., the h-overlap case and the 2h-overlap case), and independent of the overlap domain $\delta$ for the consistent overlap case (i.e., the $\delta$=0.1 and $\delta$=0.2 cases).

**Table 2. Iteration Time of Deep DDM**

| δ/h | 0.05 | 0.02 | 0.01 | 0.005 | 0.002 | 0.001 |
|---|---|---|---|---|---|---|
| h | 1382 | 1526 | 1558 | 1368 | 1923 | 2113 |
| 2h | 1227 | 1132 | 1141 | 1134 | 1623 | 2156 |
| 0.1 | 1158 | 1585 | 2170 | 1012 | 1784 | 2005 |
| 0.2 | 1833 | 1690 | 1079 | 1248 | 1488 | 1543 |

From Table 2, with the same overlap domain $\delta$, in the small overlap case (i.e., h-overlap and 2h-overlap cases), the iteration time increases with the decrease of the refinement degree h. In the consistent overlap case (i.e., $\delta$=0.1 and $\delta$=0.2 cases), the iteration time is not related to the refinement degree h.

With the same degree of refinement h, in the small overlap case (i.e., the h-overlap case and the 2h-overlap case), the iteration time is independent of the overlap domain $\delta$; in the consistent overlap case (i.e., the $\delta$=0.1 and $\delta$=0.2 cases), the iteration time is also independent of the overlap domain $\delta$.

The approximate solution of the model is as follows.

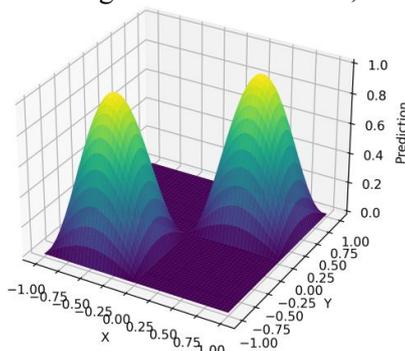

**Figure 1. Numerical Solution**

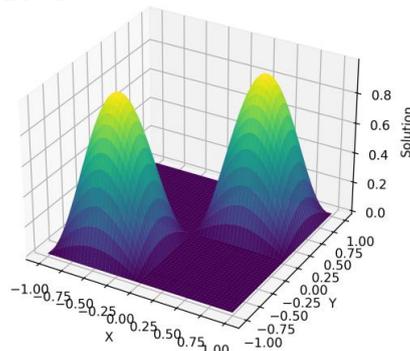

**Figure 2. Analytical Solution**

Figure 1 shows the numerical solution computed by DRPINN, Figure 2 shows the analytical solution of the problem. It can be seen that the output of the model is very close to the analytical solution of the original problem, indicating that the training of the model was successful. Next, we plot the 3D error and 2D error to check the error in detail.

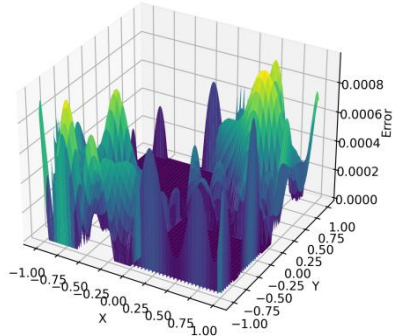

**Figure 3. Error 3d Map**

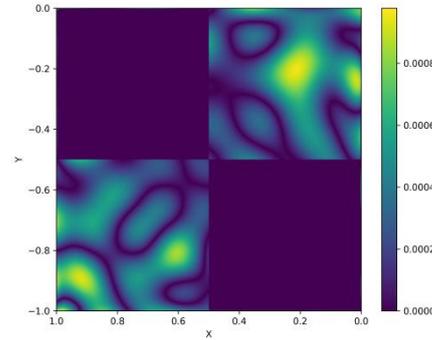

**Figure 4. Error 2D Map**

As can be seen from Figure 3 and Figure 4, although the errors show some fluctuations, this indicates that there is still room for improvement in the fitting ability of the model in some areas, and overall the errors are small, indicating that the predictive performance of the model is good.

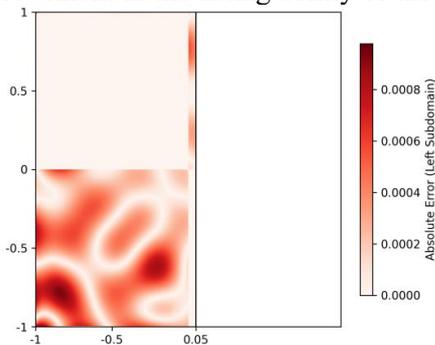
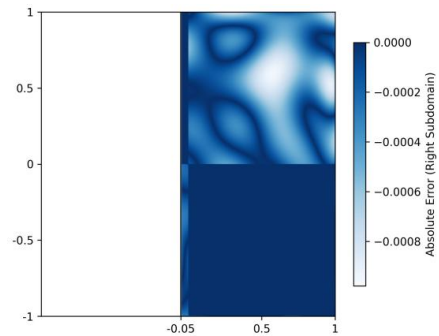

**Figure 5. Subdomain-Based Absolute Error**

Figure 5 shows the absolute error distribution between the numerical and analytical solutions, where the computational region is split using the subdomain decomposition method $x = -\delta/2$ and $x = \delta/2$, with the subdomain boundaries marked by black dashed lines. The error in the left region ($x < -\delta/2$) is in red, and the error in the right region ($x > \delta/2$) is in blue. It can be seen that the left subdomain has larger errors, mainly concentrated in the lower region, while the right subdomain has smaller and more evenly distributed errors, indicating that domain decomposition affects the error distribution and convergence characteristics.

Specific calculated data for the detailed assessment indicators are shown in Table 3.

**Table 3. Evaluation Metrics**

| MSE | MAE | REL2 | Max Error |
|---|---|---|---|
| 1.4901e-07 | 31015e-04 | 6.0813e-07 | 9.5493e-04 |

MSE stands for Mean Squared Error as followings,

$$\text{MAE} = \frac{1}{n}\sum_{i=1}^{n}|y_i - \hat{y}_i| \qquad (14)$$

MAE stands for Mean Absolute Error,

$$\text{MSE} = \frac{1}{n}\sum_{i=1}^{n}(y_i - \hat{y}_i)^2 \qquad (15)$$

Rel. L2 Error stands for Relative L2 Norm Error which is defined as,

$$\text{RelativeL2NormError} = \frac{\|y - \hat{y}\|_2}{\|y\|_2} \qquad (16)$$

Max error is defined by the infinity norm as:

$$\text{MaxError} = \max_i |y_i| \qquad (17)$$

Here, $n$ represents the number of samples, $y_i$ denotes the actual observations, and $\hat{y}_i$ represents the corresponding predicted values. Table 3 shows that the deep domain decomposition method has less error, and the MSE reaches 1.4901e-07.

## 4. Conclusion

In this paper, combining physics-informed neural network with domain decomposition method, a novel deep domain decomposition method is constructed to solve the elliptic variational inequality problem. Firstly, the computational domain is divided into multiple subdomains, then each subproblem independently trains a PINN to approximate the solution while exchanging information at interfaces. The subproblem in each sub-region is solved by using the Deep-Ritz PINN model, which the parameters in the physics-informed neural network are obtained by using adaptive moment estimation, and the areas with larger residuals can be dynamically adjusted to make the model to learn more complex regions and improving the overall training effectiveness and solution accuracy. Numerical examples show that, the new method has less mean square error and is independent of grid step length h.

**Acknowledgements**
This work is supported by the National Undergraduate Training Program for Innovation (202310595022).